\documentclass[11pt]{amsart}
\usepackage{amsmath,amssymb,amsthm,mathtools}
\usepackage{enumitem}
\usepackage{hyperref}
\usepackage[margin=1.15in]{geometry}

\newtheorem{theorem}{Theorem}

\newtheorem{lemma}{Lemma}

\title[Fejér kernels on Vilenkin groups]
{Sharp Logarithmic Bounds for Vilenkin--Fej\'er Kernels}

\author{Nika Areshidze}

\address{
Department of Mathematics,
University of California, Irvine,
CA 92697, USA
}

\email{nareshid@uci.edu}

\author{George Tephnadze}

\address{
The University of Georgia,
School of Science and Technology,
77a Merab Kostava St,
Tbilisi, 0128, Georgia
}

\email{g.tephnadze@ug.edu}



\begin{document}
\begin{abstract}
Let $G_m$ be a bounded Vilenkin group and put $R:=\sup_j m_j$. We prove the sharp uniform estimate
\[
\|K_n\|_1\leq \frac{17}{15\ln 2}\ln R,\qquad n\geq1,
\]
for the Fejér kernels associated with the Vilenkin system. The constant $17/(15\ln2)$ is best possible. In the Walsh--Paley case $R=2$, the result reduces exactly to Toledo's sharp estimate $\sup_n\|K_n\|_1=17/15$.
\end{abstract}

\maketitle

\vspace{-1em}

\section{Introduction}
The uniform boundedness of the $L^1$-norms of Fej\'er kernels plays an important role in the study of Walsh and Vilenkin systems. In 1951, Yano \cite{Yano1951} proved the classical estimate
$$\|K_n\|_1\leq 2,\quad n\geq1,$$
for the Walsh system. This bound played an important role in 
the development of the theory of norm convergence and almost-everywhere convergence for Walsh--Fej\'er means, but the exact value of the uniform bound remained unknown. In 2018, Toledo \cite{Toledo2018} determined the exact value and proved that $$ \sup_{n\geq1}\|K_n\|_1=\frac{17}{15}.$$
Thus the exact uniform $L^1$-bound was completely determined in the
Walsh case. 

More recently, Blahota \cite{Blahota2026} obtained the corresponding exact result for the Walsh--Kaczmarz system, proving that
$$
\sup_{n\geq1}\|K_n^\kappa\|_1=\frac{71}{50}.
$$
For a general fixed bounded Vilenkin system, however, the exact value of $\sup_{n\ge1}\|K_n\|_1$ remains unknown.
One of the main difficulties
is that Toledo's proof relies on positivity properties of the
Walsh--Fej\'er kernels on certain sets, whereas such positivity is not
available in the general Vilenkin setting. Consequently, Toledo's
argument does not extend directly to bounded Vilenkin systems.
 
For a bounded Vilenkin system, put
$
R:=\sup_{j\in\mathbb N}m_j<\infty.
$
It was previously shown in \cite{AnakidzeAreshidze2026} that
$$
\|K_n\|_1\leq 4R^3,\qquad n\geq1.$$  

Although this estimate provides a uniform bound, 
its dependence on $R$ is polynomial.
Moreover, in the Walsh case
$R=2$, the bound $4R^3$ is far from Toledo's sharp value $17/15$. 
In the present paper, we substantially improve the previously known polynomial estimate by establishing a logarithmic bound. Moreover, the coefficient $17/(15\ln 2)$ is optimal among all universal constants independent of the underlying bounded Vilenkin group.
In
particular, we prove the following theorem.

\begin{theorem}\label{theorem1}
Let $G_m$ be a bounded Vilenkin group and let
$$
R:=\sup_{j\geq0}m_j<\infty.
$$
Then, for every $n\geq1$,
$$
\|K_n\|_1\leq \frac{17}{15\ln2}\ln R,
$$
Moreover, the coefficient $17/(15\ln2)$ is optimal among universal constants independent of the bounded Vilenkin group.
\end{theorem}

The sharpness of the constant follows from the Walsh case. When $R=2$, Theorem~\ref{theorem1} yields
$$
\|K_n\|_1\leq \frac{17}{15},
$$
which is optimal by Toledo's result \cite{Toledo2018}. On the other
hand, a classical result of Price \cite{Price1957} implies that
$$
\sup_{n\geq1}\|K_n\|_1\gtrsim \ln R.
$$
Thus the logarithmic dependence on $R$ is optimal in order.

Our approach is different from Toledo's argument and does not rely on
positivity properties of the Fej\'er kernels. We first establish finite
cyclic estimates and an exact formula for the $L^1$-norms of the kernels
$K_{M_n}$. We then derive a one-step reduction for $\|K_n\|_1$, which,
together with an appropriate scalar inequality, allows us to complete
the proof by induction on $n$.

For the general theory of summability of Vilenkin--Fej\'er means, we
refer to the monograph of Persson, Tephnadze and Weisz \cite{PTW2022}.

The paper is organized as follows. Section~2 contains the necessary
notation and the basic identities for Dirichlet and Fej\'er kernels. In
Section~3 we establish the auxiliary finite cyclic and $L^1$ estimates
needed for the proof. Section~4 is devoted to the proof of
Theorem~\ref{theorem1}.

\section{Preliminaries and basic kernel identities}

Denote by $\mathbb{N}_{+}$ the set of positive integers, $\mathbb{N}:=%
\mathbb{N}_{+}\cup \{0\}.$ Let $m:=(m_0,m_1,\ldots)$ be a sequence of integers satisfying $m_n\ge 2,\, n\in\mathbb N.$ Denote by $Z_{m_{n}}:=\{0,1,\ldots
,m_{n}-1\}$ the additive group of integers modulo $m_{n}$. Define the group $G_m$ as the complete direct product of the groups
$Z_{m_n}$ endowed with the product topology, where each $Z_{m_n}$ has the discrete topology. Throughout the paper, we consider bounded Vilenkin groups. We put
$
R:=\sup_{n\in\mathbb{N}}m_n<\infty.
$

The direct product $\mu $ of the measures $\mu _{n}\left( \{j\}\right)
:=1/m_{n},\ (j\in Z_{m_{n}})$ is the Haar measure on $G_{m}$ with $\mu
\left( G_{m}\right) =1.$

The elements of $G_{m}$ are represented by sequences 
\begin{equation*}
x:=\left( x_{0},x_{1},\ldots,x_{n},\ldots \right),\ \left( x_{n}\in
Z_{m_{n}}\right).
\end{equation*}

A neighborhood base at $x\in G_m$ is given by
\begin{equation*}
\begin{aligned}
I_0(x)&:=G_m,\\
I_n(x)&:=\{y\in G_m:\ y_0=x_0,\ldots,y_{n-1}=x_{n-1}\},
\qquad x\in G_m,\ n\in\mathbb N.
\end{aligned}
\end{equation*}

Set $I_{n}:=I_{n}\left( 0\right) ,$ for $n\in \mathbb{N}_{+}$ \ and  $e_n=(0,\ldots,0,1,0,\ldots)$, where the $1$ occurs in the $n$-th coordinate.

\bigskip If we define the so-called generalized number system based on $m$
as follows: 
\begin{equation*}
M_{0}:=1,\ M_{n+1}:=m_{n}M_{n}\ \ \ (n\in \mathbb{N}),
\end{equation*}%
then every $n\in \mathbb{N}_{+}$ can be uniquely expressed as $%
n=\sum_{k=0}^{\infty }n_{k}M_{k},$ where $n_{k}\in Z_{m_{k}}\ (k\in \mathbb{N})$ and only finitely many of the digits $n_k$ are nonzero. Let 
$
\left\vert n\right\vert :=\max \{j\in \mathbb{N}:n_{j}\neq 0\},
$
that is $M_{\left\vert n\right\vert }\leq n< M_{\left\vert n\right\vert
+1}.$ 

Next, we introduce the Vilenkin system on $G_m$. First, we define the generalized Rademacher functions
$r_k:G_m\rightarrow\mathbb C$ by
\begin{equation*}
r_{k}\left( x\right) :=\exp \left( 2\pi ix_{k}/m_{k}\right) ,\text{ }\left(
i^{2}=-1,x\in G_{m},\text{ }k\in \mathbb{N}\right) .
\end{equation*}

Now define the Vilenkin system $\,\,\,\psi :=(\psi _{n}:n\in\mathbb{N})$ on 
$G_{m}$ by
\begin{equation*}
\psi _{n}(x):=\prod\limits_{k=0}^{\infty }r_{k}^{n_{k}}\left( x\right) \ \ \
\left( n\in\mathbb{N}\right).
\end{equation*}

When $m\equiv 2$, the Vilenkin system reduces to the Walsh--Paley system. The Vilenkin system is orthonormal and complete in $L_{2}\left( G_{m}\right) 
$ (see, e.g., \cite{AVD}).

If $f\in L_{1}\left( G_{m}\right) $ we can define Fourier coefficients,
partial sums, Dirichlet kernels, Fejér means, and Fejér kernels with respect to
the Vilenkin system in the usual manner: 
\begin{equation*}
\widehat{f}\left( k\right) :=\int_{G_{m}}f\overline{\psi }_{k}d\mu \text{%
\thinspace\ \ \ }\left( \text{ }k\in \mathbb{N}\right) ,
\end{equation*}%
\begin{eqnarray*}
S_{n}f &:&=\sum_{k=0}^{n-1}\widehat{f}\left( k\right) \psi _{k},\ \text{%
\qquad }D_{n}:=\sum_{k=0}^{n-1}\psi _{k\text{ }}\text{ \ \ \ \ \ \ }\left( 
\text{ }n\in \mathbb{N}_{+}\text{ }\right) , \\
\sigma _{n}f &:&=\frac{1}{n}\sum_{k=1}^{n}S_{k}f,\text{\qquad\ \ \ }K_{n}:=%
\frac{1}{n}\sum_{k=1}^{n}D_{k}\text{ \ \ }\left( \text{
}n\in \mathbb{N}_{+}\text{ }\right) .
\end{eqnarray*}

Recall that (see e.g., \cite{AVD}) 
\begin{equation}
\quad \hspace*{0in}D_{M_{n}}\left( x\right) =\left\{ 
\begin{array}{l}
\text{ }M_{n},\text{ \thinspace \thinspace if \thinspace \thinspace }x\in
I_{n}, \\ 
\text{ }0,\text{ \ \ \thinspace \thinspace \thinspace if \ \thinspace
\thinspace }x\notin I_{n},%
\end{array}%
\right.  \label{3}
\end{equation}%
and \ 
\begin{equation}
D_{s_{n}M_{n}}=D_{M_{n}}\sum_{k=0}^{s_{n}-1}\psi
_{kM_{n}}=D_{M_{n}}\sum_{k=0}^{s_{n}-1}r_{n}^{k},  \label{9dn}
\end{equation}%
where $n\in \mathbb{N}$ and $1\leq s_{n}\leq m_{n}-1.$

The following three identities for Fej\'er kernels can be found in \cite{PTW2022}.
Let $n,t\in\mathbb N$ with $n>t$ and $x\in I_t\setminus I_{t+1}$. Then
\begin{equation}\label{lemma3}
K_{M_{n}}\left( x\right) =\left\{ 
\begin{array}{ll}
0, & \text{if }x-x_{t}e_{t}\notin I_{n}, \\ 
\frac{M_{t}}{1-r_{t}\left( x\right) }, & \text{if }x-x_{t}e_{t}\in I_{n}.%
\end{array}%
\right.
\end{equation}

Let $n\in \mathbb{N}$ and $1\leq s\leq m_n-1$. Then
\begin{equation} \label{lemma6kn}
sM_nK_{sM_n}=\sum_{l=0}^{s-1}\left(\sum_{i=0}^{l-1}r_n^i\right)M_nD_{M_n}+\left(\sum_{l=0}^{s-1}r_n^{l}\right)M_nK_{M_n}
\end{equation}
Let
$
n=sM_k+q,
\quad 1\le s<m_k,
\quad 0\le q<M_k.
$
Then
\begin{equation}\label{eq:one-digit-fejer}
nK_n=sM_kK_{sM_k}+qD_{sM_k}+r_k^sqK_q.
\end{equation}

\section{Auxiliary $L^1$ estimates}

For $m\ge2$ and $1\le t\le m$, put
\[
\zeta_m:=e^{2\pi i/m},
\qquad
B_t(z):=\sum_{j=0}^{t-1}z^j,
\]
and, for $1\le s<m$,
\[
A_s(z):=\sum_{\ell=1}^{s-1}B_\ell(z).
\]
We shall use the quantities
\begin{equation}\label{eq:lambda-def}
\lambda_m(t):=\frac1m\sum_{a=0}^{m-1}|B_t(\zeta_m^a)|,
\end{equation}
\begin{equation}\label{eq:alpha-def}
\alpha_m(s):=\frac1{sm}\sum_{a=0}^{m-1}|A_s(\zeta_m^a)|,
\end{equation}

For $m\ge2$, put
\begin{equation}\label{eq:beta-def}
\beta_m:=\frac1m\sum_{a=1}^{m-1}\frac1{|1-\zeta_m^a|}
=\frac1{2m}\sum_{a=1}^{m-1}\csc\frac{\pi a}{m},
\end{equation}
and, for $R\ge2$, set
\begin{equation}\label{eq:PR-def}
P_R:=\max_{2\le m\le R}\frac{m}{m-1}\beta_m.
\end{equation}

\begin{lemma}[Kouba's cosecant estimate \cite{Kouba2016}]\label{lemma:kouba}
Let $\gamma$ denote the Euler--Mascheroni constant,
	\[
	\gamma:=\lim_{n\to\infty}
	\left(
	\sum_{k=1}^{n}\frac{1}{k}-\ln n
	\right), \qquad (\gamma \approx 0.5772156649)
\]
and let
	\[
	\eta:=\gamma-\ln\frac\pi2,\qquad (	
	\eta \approx0.1256329596.
	)\]	
	
	$$
b(x):=\frac{\ln x+\eta}{\pi},
	\qquad x>1,$$
	
	\[
	\Pi(x):=\frac{x}{x-1}b(x).
	\]
	Then, for every integer $m\ge2$,
	\begin{equation}\label{eq:beta-kouba}
	\beta_m<b(m).
	\end{equation}
	Moreover, $\Pi$ is increasing on $[2,\infty)$, and
	\begin{equation}\label{eq:Pi-derivative-upper}
	\Pi'(x)<\frac1{\pi x},
	\qquad x\ge3.
	\end{equation}
	Consequently,
	\begin{equation}\label{eq:PR-Pi}
	P_R<\Pi(R).
	\end{equation}
\end{lemma}

\begin{proof}
	Kouba's cosecant inequality \cite[Theorem~3.4, with $n=0$]{Kouba2016}
	states that
	\[
	\sum_{a=1}^{m-1}\csc\frac{\pi a}{m}
	<\frac{2m}{\pi}
	\left(\ln m+\gamma-\ln\frac\pi2\right).
	\]
	Together with \eqref{eq:beta-def}, this gives \eqref{eq:beta-kouba}.
	
	A direct differentiation gives
	\begin{equation}\label{eq:Pi-derivative}
	\Pi'(x)
	=\frac{x-1-\ln x-\eta}{\pi(x-1)^2}.
	\end{equation}
	The numerator is increasing for $x>1$ and is positive at $x=2$; hence
	$\Pi$ is increasing on $[2,\infty)$.  Furthermore,
	\eqref{eq:Pi-derivative-upper} is equivalent to
	\[
	\ln x+\eta>1-\frac1x,
	\]
	which holds for $x\ge3$.  Finally, for every integer $2\le m\le R$,
	\[
	\frac{m}{m-1}\beta_m
	<\frac{m}{m-1}b(m)=\Pi(m)\le\Pi(R),
	\]
	and \eqref{eq:PR-Pi} follows.
\end{proof}

\begin{lemma}[Finite cyclic estimates]\label{lemma:cyclic-estimates}
	Let $m\ge2$. Then, for $1\le t\le m$,
	\begin{equation}\label{eq:lambda-bounds}
	1\le\lambda_m(t)\le\sqrt t,
	\end{equation}
	and
	\begin{equation}\label{eq:lambda-beta}
	\lambda_m(t)\le\frac tm+2\beta_m.
	\end{equation}
	Moreover,
	\begin{equation}\label{eq:lambda2}
	\lambda_m(2)\le\frac43,
	\end{equation}
	while, for $1\le s<m$,
	\begin{equation}\label{eq:alpha-beta}
	\alpha_m(s)
	\le\frac{s-1}{2m}+2\frac{s-1}{s}\beta_m.
	\end{equation}
	In particular, for $m\ge3$,
	\begin{equation}\label{eq:alpha2}
	\alpha_m(2)=\frac12.
	\end{equation}
\end{lemma}

\begin{proof}
Since
\[
\begin{aligned}
\left|B_t(\zeta_m^a)\right|^2
&=
\left(\sum_{j=0}^{t-1}\zeta_m^{aj}\right)
\left(\sum_{\ell=0}^{t-1}\zeta_m^{-a\ell}\right)=
\sum_{j=0}^{t-1}\sum_{\ell=0}^{t-1}
\zeta_m^{a(j-\ell)},
\end{aligned}
\]
averaging over $a=0,\ldots,m-1$, we obtain
\[
\frac1m\sum_{a=0}^{m-1}
\left|B_t(\zeta_m^a)\right|^2
=
\sum_{j=0}^{t-1}\sum_{\ell=0}^{t-1}
\left(
\frac1m\sum_{a=0}^{m-1}
\zeta_m^{a(j-\ell)}
\right).
\]

Now use the standard orthogonality relation for the $m$-th roots of unity:
\[
\frac1m\sum_{a=0}^{m-1}\zeta_m^{ar}
=
\begin{cases}
1, & r\equiv 0 \pmod m,\\[1mm]
0, & r\not\equiv 0 \pmod m.
\end{cases}
\]	

Since $1\le t\le m$,
	$
	|j-\ell|\le t-1<m,
$
	we obtain
	$
	j-\ell\equiv0\pmod m
	\quad\Longleftrightarrow\quad
	j=\ell.
	$
	Thus
	$$
	\begin{aligned}
	\frac1m\sum_{a=0}^{m-1}|B_t(\zeta_m^a)|^2=
	\sum_{j=0}^{t-1}\sum_{\ell=0}^{t-1}\delta_{j\ell}=
	\sum_{j=0}^{t-1}1=t.
	\end{aligned}
	$$
The Cauchy--Schwarz inequality gives
\[
\begin{aligned}
\lambda_m(t)
&=
\frac1m\sum_{a=0}^{m-1}
\left|B_t(\zeta_m^a)\right|\le
\frac1m
\left(
\sum_{a=0}^{m-1}1^2
\right)^{1/2}
\left(
\sum_{a=0}^{m-1}
\left|B_t(\zeta_m^a)\right|^2
\right)^{1/2}\\
&=
\frac1m\sqrt{m}\sqrt{mt}=
\sqrt{t}.
\end{aligned}
\]
Therefore,
$
	\lambda_m(t)\le \sqrt{t}.
$
	On the other hand,
	$$
	\frac1m\sum_{a=0}^{m-1}B_t(\zeta_m^a)=1,
	$$
	and therefore $\lambda_m(t)\ge1$.

Among the $m$-th roots of unity
$\zeta_m^a,$ 
$ a=0,\ldots,m-1,$
the case $a=0$ is special because
$
\zeta_m^0=1.
$
Thus, the term corresponding to $a=0$ is
\[
|B_t(\zeta_m^0)|=|B_t(1)|=t.
\]
Hence, separating the root $1$ means taking the term corresponding to
$a=0$ out of the sum:
\[
\begin{aligned}
\lambda_m(t)=
\frac1m|B_t(1)|
+\frac1m\sum_{a=1}^{m-1}|B_t(\zeta_m^a)|=
\frac{t}{m}
+\frac1m\sum_{a=1}^{m-1}|B_t(\zeta_m^a)|.
\end{aligned}
\]

For $a\neq0$, we have $\zeta_m^a\neq1$, and therefore the geometric-sum
formula gives
\begin{eqnarray*}
B_t(\zeta_m^a)
=
\frac{1-\zeta_m^{at}}{1-\zeta_m^a}.
\end{eqnarray*}
Hence, using the estimate
\begin{eqnarray*}
\left|\frac{1-z^t}{1-z}\right|\le\frac2{|1-z|}
\end{eqnarray*}
we obtain
\[
	\lambda_m(t)
	=
	\frac{t}{m}
	+
	\frac1m\sum_{a=1}^{m-1}
	\frac{|1-\zeta_m^{at}|}{|1-\zeta_m^a|}\le\frac tm+2\beta_m,
\]
which gives
\eqref{eq:lambda-beta}.
	
	To prove \eqref{eq:lambda2}, note that for $t=2$, $ B_2(z)=1+z. $
Hence, for \(z=\zeta_m^a=e^{2\pi ia/m}\),
$ B_2(\zeta_m^a) = 1+e^{2\pi ia/m}. $
Now use
$$ 1+e^{i\theta} = e^{i\theta/2} \left( e^{-i\theta/2}+e^{i\theta/2} \right) = 2e^{i\theta/2}\cos\frac{\theta}{2}. $$
Taking absolute values gives
$ |1+e^{i\theta}| = 2\left|\cos\frac{\theta}{2}\right|, $
because
$ |e^{i\theta/2}|=1. $
With $ \theta=\frac{2\pi a}{m}, $
we obtain
$ |B_2(\zeta_m^a)| = \left|1+e^{2\pi ia/m}\right| = 2\left|\cos\frac{\pi a}{m}\right|. $

Therefore, by definition,
	\[
	\lambda_m(2)
	=\frac2m\sum_{a=0}^{m-1}\left|\cos\frac{\pi a}{m}\right|.
	\]
	
We consider separately the cases where $m$ is even and where $m$ is odd. In both cases we use the finite cosine-sum identity (Dirichlet kernel identity)
\begin{eqnarray}\label{2}
1+2\sum_{a=1}^{n}\cos(a\theta)
=
\frac{\sin\left(\left(n+\frac12\right)\theta\right)}
{\sin\left(\frac{\theta}{2}\right)}
\end{eqnarray}
and the identity
\begin{eqnarray}\label{1}
\cos\frac{\pi(m-a)}m = -\cos\frac{\pi a}{m}, 
\end{eqnarray}

Suppose that \(m\) is even. By \eqref{1} the absolute values are symmetric. 
Moreover,
$
\cos\frac{\pi(m/2)}{m}
=
\cos\frac{\pi}{2}
=
0.
$
By using \eqref{2} we obtain
\[
\begin{aligned}
1+2\sum_{a=1}^{m/2-1}\cos\frac{\pi a}{m}=
\frac{
	\sin\left(\left(\frac{m}{2}-\frac12\right)\frac{\pi}{m}\right)
}{
	\sin\frac{\pi}{2m}
}=
\frac{
	\sin\left(\frac{\pi}{2}-\frac{\pi}{2m}\right)
}{
	\sin\frac{\pi}{2m}
}=
\frac{
	\cos\frac{\pi}{2m}
}{
	\sin\frac{\pi}{2m}
}=
\cot\frac{\pi}{2m}.
\end{aligned}
\]
Hence,
\[
\lambda_m(2)
=
\frac{2}{m}\cot\frac{\pi}{2m}.
\]
	
	To estimate this quantity, we use
	$
	\tan x>x,$ $ 0<x<\frac{\pi}{2},
	$
	which is equivalent to
	$
	\cot x<\frac1x.
	$
	Thus, taking \(x=\pi/(2m)\), we obtain
	\[
	\cot\frac{\pi}{2m}
	<
	\frac{2m}{\pi}.
	\]
	Consequently,
	\[
	\lambda_m(2)
	=
	\frac2m\cot\frac{\pi}{2m}
	<
	\frac2m\frac{2m}{\pi}
	=
	\frac4\pi<\frac43.
	\]
If $m$ is odd, then combining \eqref{2} and \eqref{1} gives
	\[
	\begin{aligned}
	1+2\sum_{a=1}^{(m-1)/2}\cos\frac{\pi a}{m}=
	\frac{
		\sin\left(
		\left(\frac{m-1}{2}+\frac12\right)\frac{\pi}{m}
		\right)
	}{
		\sin\frac{\pi}{2m}
	}=
	\frac{
		\sin\frac{\pi}{2}
	}{
		\sin\frac{\pi}{2m}
	}=
	\frac{1}{
		\sin\frac{\pi}{2m}
	}.
	\end{aligned}
	\]
	
	Consequently,	
	\[
	\lambda_m(2)=\frac{2}{m\sin(\pi/(2m))}.
	\]
	For $m\ge3$, concavity of $\sin x$ on $[0,\pi/6]$ gives
	$\sin x\ge 3x/\pi$, and hence, 
	\[
	\sin\frac{\pi}{2m}
	\ge
	\frac{3}{\pi}\frac{\pi}{2m}
	=
	\frac{3}{2m}.
	\]
	Consequently,
	\[
	\lambda_m(2)
	\le
	\frac{2}{m\cdot\frac{3}{2m}}
	=
	\frac{4}{3}.
	\]

	Finally, since $A_s=\sum_{\ell=1}^{s-1}B_\ell$, substituting the representation of \(A_s\) in \eqref{eq:alpha-def}, we get
	\begin{eqnarray}
	\alpha_m(s) &=& \frac{1}{sm} \sum_{a=0}^{m-1} \left| \sum_{\ell=1}^{s-1} B_\ell(\zeta_m^a) \right|\le \frac{1}{sm} \sum_{a=0}^{m-1} \sum_{\ell=1}^{s-1} \left|B_\ell(\zeta_m^a)\right|\\ \notag
	&=& \frac1s \sum_{\ell=1}^{s-1} \left( \frac1m \sum_{a=0}^{m-1} \left|B_\ell(\zeta_m^a)\right| \right). 
	\end{eqnarray}

	By definition,	
	$\lambda_m(\ell) = \frac1m \sum_{a=0}^{m-1} |B_\ell(\zeta_m^a)|, $	
	hence
	
	$$  \alpha_m(s) \le \frac1s\sum_{\ell=1}^{s-1}\lambda_m(\ell). $$

	Applying \eqref{eq:lambda-beta}, we obtain
	\[
	\begin{aligned}
	\alpha_m(s)\le\frac1s\sum_{\ell=1}^{s-1}\lambda_m(\ell)\le\frac1s\sum_{\ell=1}^{s-1}
	\left(\frac\ell m+2\beta_m\right)
	=\frac{s-1}{2m}+2\frac{s-1}{s}\beta_m.
	\end{aligned}
	\]
	Since $A_2\equiv1$, \eqref{eq:alpha2} follows.
\end{proof}

\begin{lemma}\label{lemma7kn}
	Let $n\geq 0$. Then
	\begin{equation}\label{eq:KMn-L1-exact}
	\|K_{M_n}\|_1
	=\frac{M_n+1}{2M_n}
	+\sum_{t=0}^{n-1}\frac{M_{t+1}}{M_n}\beta_{m_t},
	\end{equation}
	and consequently
	\begin{equation}\label{eq:KMn-L1-bound}
	\|K_{M_n}\|_1
	\leq\frac{M_n+1}{2M_n}
	+P_R\left(1-\frac1{M_n}\right)
	\leq P_R+\frac12.
	\end{equation}
\end{lemma}

\begin{proof}
By \eqref{lemma3}, $K_{M_n}$ vanishes outside
\[
\mathcal S_n:=I_n\cup
\bigcup_{t=0}^{n-1}\bigcup_{a=1}^{m_t-1}I_n(ae_t),
\]
while for $x\in I_n(ae_t)$, $0\leq t<n$ and $1\leq a<m_t$,
\[
K_{M_n}(x)=\frac{M_t}{1-r_t(x)}=\frac{M_t}{1-\zeta_{m_t}^a}.
\]
On $I_n$ one has $D_j(x)=j$ for $1\leq j\leq M_n$, and hence
\[
K_{M_n}(x)=\frac1{M_n}\sum_{j=1}^{M_n}D_j(x)=\frac{M_n+1}{2}.
\]
Since $\mu(I_n)=\mu(I_n(ae_t))=1/M_n$, we obtain the integral
\[
\begin{aligned}
\|K_{M_n}\|_1
&=\frac{M_n+1}{2M_n}
+\sum_{t=0}^{n-1}\frac1{M_n}\sum_{a=1}^{m_t-1}
\frac{M_t}{|1-\zeta_{m_t}^a|}\\
&=\frac{M_n+1}{2M_n}
+\sum_{t=0}^{n-1}\frac{M_{t+1}}{M_n}\beta_{m_t},
\end{aligned}
\]
which proves \eqref{eq:KMn-L1-exact}.

We now estimate this $L^1$-norm. By the definition of $P_R$,
\[
\beta_{m_t}\leq\frac{m_t-1}{m_t}P_R.
\]
Therefore
\[
\begin{aligned}
\sum_{t=0}^{n-1}\frac{M_{t+1}}{M_n}\beta_{m_t}
&\leq \frac{P_R}{M_n}\sum_{t=0}^{n-1}M_t(m_t-1)\\
&=\frac{P_R}{M_n}\sum_{t=0}^{n-1}(M_{t+1}-M_t)=P_R\left(1-\frac1{M_n}\right).
\end{aligned}
\]
Hence
\[
\|K_{M_n}\|_1
\leq\frac{M_n+1}{2M_n}
+P_R\left(1-\frac1{M_n}\right).
\]
Finally, since
\[
P_R\geq\frac{2}{2-1}\beta_2=\frac12,
\]
we have
\[
\frac{M_n+1}{2M_n}+P_R\left(1-\frac1{M_n}\right)
=P_R+\frac12+\frac{1/2-P_R}{M_n}
\leq P_R+\frac12,
\]
which proves \eqref{eq:KMn-L1-bound}.
\end{proof}

\begin{lemma}[One-step $L^1$ reduction]\label{lemma:L1-reduction}
Let
$
n=sM_k+q,
\ 1\le s<m_k,
\ 0\le q<M_k.
$
Then
\begin{equation}\label{eq:induction-reduction}
\begin{aligned}
n\|K_n\|_1\le q\|K_q\|_1+sM_k\Bigg[\alpha_{m_k}(s)+(1+P_R)\frac{\lambda_{m_k}(s)}s
+\frac{\lambda_{m_k}(s+1)-1}{2s}\Bigg].
\end{aligned}
\end{equation}
When $q=0$, all terms containing $qK_q$ are understood to be zero; in this case $b=v=0$.
\end{lemma}

\begin{proof}
By combining \eqref{9dn}, \eqref{lemma6kn} and \eqref{eq:one-digit-fejer} we obtain
$
D_{sM_k}=D_{M_k}B_s(r_k)
$
and
\[
sM_kK_{sM_k}
=M_kD_{M_k}A_s(r_k)+M_kK_{M_k}B_s(r_k).
\]
Consequently,
\begin{equation}\label{eq:master-decomposition}
nK_n
=M_kD_{M_k}A_s(r_k)
+(M_kK_{M_k}+qD_{M_k})B_s(r_k)
+r_k^sqK_q.
\end{equation}

Write
\[
U:=M_kD_{M_k}A_s(r_k)
+(M_kK_{M_k}+qD_{M_k})B_s(r_k).
\]
By \eqref{lemma3}, $K_{M_k}$ is supported in $\mathcal S_k$.
For $k=0$ this is immediate, since $K_{M_0}=K_1\equiv1$ and $\mathcal S_0=G_m$.
On $\mathcal S_k\setminus I_k$ we have $D_{M_k}=0$, and hence
\[
U=M_kK_{M_k}B_s(r_k).
\]

Since $K_{M_k}$ depends only on the coordinates
$x_0,\ldots,x_{k-1}$,
whereas $r_k(x)$, and therefore $B_s(r_k(x))$, depends only on the
$k$-th coordinate $x_k$, these two factors can be separated with
respect to the product Haar measure. Consequently, by the definition
of $\lambda_{m_k}(s)$,
\[
\begin{aligned}
\int_{\mathcal S_k\setminus I_k}|U|\,d\mu
&=
M_k
\left(
\int_{\mathcal S_k\setminus I_k}
|K_{M_k}(x)|\,d\mu(x)
\right)
\left(
\frac1{m_k}\sum_{a=0}^{m_k-1}
|B_s(\zeta_{m_k}^{a})|
\right)\\
&=M_k\lambda_{m_k}(s)
\int_{\mathcal S_k\setminus I_k}
|K_{M_k}(x)|\,d\mu(x).
\end{aligned}
\]

Next, \(K_{M_k}\) is supported on \(\mathcal S_k\). Therefore,
\[
\int_{\mathcal S_k\setminus I_k}|K_{M_k}|\,d\mu
=
\|K_{M_k}\|_1
-
\int_{I_k}|K_{M_k}|\,d\mu.
\]

On \(I_k\),
\[
K_{M_k}(x)=\frac{M_k+1}{2}, \ \ \ \text{
and} \ \ \
\mu(I_k)=\frac1{M_k}.
\]
Hence,
\[
\int_{I_k}|K_{M_k}|\,d\mu
=
\frac{M_k+1}{2}\frac1{M_k}
=
\frac{M_k+1}{2M_k}.
\]

It follows that
\[
\int_{\mathcal S_k\setminus I_k}|K_{M_k}|\,d\mu
=
\|K_{M_k}\|_1
-
\frac{M_k+1}{2M_k},
\]
and therefore
\begin{equation}\label{eq:outside-U-fullnorm}
\begin{aligned}
\int_{\mathcal S_k\setminus I_k}|U|\,d\mu
&=M_k\lambda_{m_k}(s)
\left(\|K_{M_k}\|_1-\frac{M_k+1}{2M_k}\right)\\
&\le (M_k-1)P_R\lambda_{m_k}(s).
\end{aligned}
\end{equation}

We next exploit the exact cancellation on $I_k$. 
If $1\le q<M_k$, then
$D_j=j$ on $I_k$ for $1\le j\le q$, whence
\[
K_q=\frac{q+1}{2}\qquad\text{on }I_k.
\]
Set
\[
H:=\frac{M_k(M_k+1)}2+qM_k,
\qquad
b:=\frac{q(q+1)}2,
\qquad
h:=H-b.
\]
Since $0\le q\le M_k-1$,
\begin{equation}\label{eq:h-b-ranges}
0\le h\le M_k^2,
\qquad
0\le b\le\frac{M_k^2}{2}.
\end{equation}
Indeed, if we define 
\[
h(q)
:=
\frac{M_k(M_k+1)}{2}
+qM_k
-\frac{q(q+1)}{2}, \ \ \text{
where} \ \
0\le q\le M_k-1,
\]
then
\[
\begin{aligned}
h(q+1)-h(q)&=M_k-\frac{(q+1)(q+2)-q(q+1)}{2}=M_k-\frac{2q+2}{2}\\
&=M_k-q-1\geq 0,
\end{aligned}
\]
while $h(0)=M_k(M_k+1)/2$ and $h(M_k-1)=M_k^2$.

On $I_k$, formula \eqref{eq:master-decomposition} becomes
\[
M_k^2A_s(r_k)+HB_s(r_k)+br_k^s.
\]
Using $H=h+b$ and $B_s(z)+z^s=B_{s+1}(z)$, we obtain
\begin{equation}\label{eq:central-identity}
M_k^2A_s(z)+HB_s(z)+bz^s
=M_k^2A_s(z)+hB_s(z)+bB_{s+1}(z).
\end{equation}
Put
\begin{equation}\label{eq1}
u:=\frac{h}{M_k^2}\in[0,1],
\qquad
v:=\frac{b}{M_k^2}\in\left[0,\frac12\right].
\end{equation}

Since $U$ is supported in $\mathcal S_k$, the triangle inequality outside
$I_k$, together with \eqref{eq:outside-U-fullnorm}, gives
\[
\begin{aligned}
n\|K_n\|_1
&\le q\|K_q\|_1+(M_k-1)P_R\lambda_{m_k}(s)\\
&\quad+\int_{I_k}|nK_n|\,d\mu-q\int_{I_k}|K_q|\,d\mu.
\end{aligned}
\]
Furthermore,
\[
q\int_{I_k}|K_q|\,d\mu=M_kv,
\]
and \eqref{eq:central-identity} yields
\[
\int_{I_k}|nK_n|\,d\mu
=\frac{M_k}{m_k}\sum_{a=0}^{m_k-1}
\left|
A_s(\zeta_{m_k}^a)
+uB_s(\zeta_{m_k}^a)
+vB_{s+1}(\zeta_{m_k}^a)
\right|.
\]
Hence
\begin{eqnarray}\label{eq:bookkeeping-fullnorm}
n\|K_n\|_1
&\le& q\|K_q\|_1+(M_k-1)P_R\lambda_{m_k}(s)\notag\\
&+&M_k\left[
\frac1{m_k}\sum_{a=0}^{m_k-1}
\left|
A_s(\zeta_{m_k}^a)
+uB_s(\zeta_{m_k}^a)
+vB_{s+1}(\zeta_{m_k}^a)
\right|-v
\right].
\end{eqnarray}

Using the standard orthogonality relation for the $m_k$-th roots of unity,
\[
\frac{1}{m_k}\sum_{a=0}^{m_k-1}\zeta_{m_k}^{aj}
=
\begin{cases}
1, & j\equiv 0 \pmod{m_k},\\[1mm]
0, & j\not\equiv 0 \pmod{m_k}.
\end{cases}
\]
Since
$
0\le j\le s\le m_k-1,
$
the only value of $j$ divisible by $m_k$ is $j=0$; therefore,
\[
\begin{aligned}
\frac{1}{m_k}\sum_{a=0}^{m_k-1}
B_{s+1}(\zeta_{m_k}^a)
&=
\frac{1}{m_k}\sum_{a=0}^{m_k-1}
\sum_{j=0}^{s}\zeta_{m_k}^{aj}=
\sum_{j=0}^{s}
\left(
\frac{1}{m_k}\sum_{a=0}^{m_k-1}
\zeta_{m_k}^{aj}
\right)=1.
\end{aligned}
\]
Hence
$$\lambda_{m_k}(s+1)=
\frac{1}{m_k}\sum_{a=0}^{m_k-1}
\left|B_{s+1}(\zeta_{m_k}^a)\right|\geq \left|\frac{1}{m_k}\sum_{a=0}^{m_k-1}
B_{s+1}(\zeta_{m_k}^a)\right|=1.$$ 
Therefore, using \eqref{eq1} and the triangle inequality,
\[
\begin{aligned}
&\frac1{m_k}\sum_{a=0}^{m_k-1}
\left|
A_s(\zeta_{m_k}^a)
+uB_s(\zeta_{m_k}^a)
+vB_{s+1}(\zeta_{m_k}^a)
\right|-v\\
&\qquad\le
s\alpha_{m_k}(s)+\lambda_{m_k}(s)
+\frac12\bigl(\lambda_{m_k}(s+1)-1\bigr).
\end{aligned}
\]
Substituting this into \eqref{eq:bookkeeping-fullnorm}, we obtain
\[
\begin{aligned}
n\|K_n\|_1
&\le q\|K_q\|_1+sM_k\Bigg[\alpha_{m_k}(s)\\
&\qquad+\left(1+P_R\left(1-\frac1{M_k}\right)\right)
\frac{\lambda_{m_k}(s)}s
+\frac{\lambda_{m_k}(s+1)-1}{2s}\Bigg].
\end{aligned}
\]
Since $1-1/M_k\le1$, \eqref{eq:induction-reduction} follows.
\end{proof}

\begin{lemma}[The scalar inequality]\label{lemma:scalar}
	Let $R\ge3$, $2\le m\le R$, and $1\le s<m$. Then
	\begin{equation}\label{eq:master-inequality}
	\alpha_m(s)
	+(1+P_R)\frac{\lambda_m(s)}s
	+\frac{\lambda_m(s+1)-1}{2s}
	\le \frac{17}{15\ln2}\ln R.
	\end{equation}
\end{lemma}

\begin{proof}
	We consider three cases, where we need the following estimates:
	$$0.6931<\ln2<0.6932, \ 1.0986<\ln3<1.0987, \ 3.1415<\pi<3.1416,$$
	and 
	$$0.1256<\eta<0.1257.$$
	
	\smallskip
	\noindent
	\emph{Case 1: $s=1$.}
	Here $\alpha_m(1)=0$ and $\lambda_m(1)=1$. By
	\eqref{eq:lambda2} and \eqref{eq:PR-Pi},
	$$
	\alpha_m(1)+(1+P_R)\lambda_m(1)
	+\frac{\lambda_m(2)-1}{2}
	\le \frac76+\Pi(R).
	$$
	Set
	$$
	F_1(R):=\frac{17}{15\ln2}\ln R-\frac76-\Pi(R).
	$$
	By \eqref{eq:Pi-derivative-upper},
	$$
	F_1'(R)>
	\frac1R\left(\frac{17}{15\ln2}-\frac1\pi\right)>0,
	\qquad R\ge3.
	$$
	Thus $F_1$ is increasing on $[3,\infty)$. Moreover,
	$$
	F_1(3)
	=
	\frac{17\ln3}{15\ln2}
	-\frac76
	-\frac32\,\frac{\ln3+\eta}{\pi},
	\qquad
	\eta=\gamma-\ln\frac{\pi}{2}.
	$$

	Since
	\begin{equation}\label{eq5}
	\frac{17\ln3}{15\ln2}\approx1.796,
	\qquad
	\Pi(3)=
	\frac32\,b(3)=\frac32\,\frac{\ln3+\eta}{\pi}\approx0.585,
    \end{equation}
	we obtain
	$$
	F_1(3)\approx1.796-1.167-0.585\approx0.045
	$$
	and
	$$F_1(3)>0.044>0.$$
	Therefore,
	$
	F_1(R)\ge F_1(3)>0,
	\ R\ge3.
    $

%%%%%%%%%%%%%%%%%%%%%%%%%%%%%%%%%%%%%%%%%%%%%%%%%%%%%%%%%%%%%%%%%%%%%%%%%%%%%%%%%%%%%%%
	
	\smallskip
	\noindent
	\emph{Case 2: $s=2$.}
	
Now $m\ge3$. By \eqref{eq:alpha2}, \eqref{eq:lambda2},
\eqref{eq:lambda-bounds}, and \eqref{eq:PR-Pi},
\[
\begin{aligned}
&\alpha_m(2)
+(1+P_R)\frac{\lambda_m(2)}2
+\frac{\lambda_m(3)-1}{4}\\
&\qquad\le
\frac12+\frac23(1+\Pi(R))
+\frac{\sqrt3-1}{4}.
\end{aligned}
\]
Set
\[
F_2(R):=
\frac{17}{15\ln2}\ln R
-\frac12
-\frac23(1+\Pi(R))
-\frac{\sqrt3-1}{4}.
\]
By \eqref{eq:Pi-derivative-upper},
\[
F_2'(R)
>
\frac1R
\left(
\frac{17}{15\ln2}
-\frac{2}{3\pi}
\right)
>0,
\qquad R\ge3.
\]
Thus $F_2$ is increasing on $[3,\infty)$. At $R=3$,
\[
\begin{aligned}
F_2(3)
&=\frac{17\ln3}{15\ln2}
-\frac{\ln3+\eta}{\pi}
-\frac{11+3\sqrt3}{12}\\
&=
\frac{17\ln3}{15\ln2}
-\frac12
-\frac23(1+\Pi(3))
-\frac{\sqrt3-1}{4}.
\end{aligned}
\]

By combining  \eqref{eq5} and
\[
\frac{\sqrt3-1}{4}\approx0.183,
\]
we obtain
\[
\begin{aligned}
F_2(3)
&\approx
1.796
-0.5
-\frac23(1+0.585)
-0.183\approx0.057
\end{aligned}
\]
and
$$F_2(3)>0.056>0.$$
Consequently,
$
F_2(R)\ge F_2(3)>0,
\ R\ge3,
$
which proves \eqref{eq:master-inequality} for $s=2$.

%%%%%%%%%%%%%%%%%%%%%%%%%%%%%%%%%%%%%%%%%%%%%%%%%%%%%%%%%%%%%%%%%%%%%%%%%%%%%%%%%%%%

\smallskip
\noindent
\emph{Case 3: $s\ge3$.}
Then $R\ge4$. By \eqref{eq:alpha-beta}, \eqref{eq:beta-kouba}, and the
monotonicity of $b$,
\begin{equation}\label{eq:alpha-R}
\alpha_m(s)
\le
\frac{R-2}{2R}
+2\frac{R-2}{R-1}b(R).
\end{equation}
Moreover, by \eqref{eq:lambda-bounds},
\[
\frac{\lambda_m(s)}s\le\frac1{\sqrt3},
\]
and
\[
\frac{\lambda_m(s+1)-1}{2s}
\le
\frac1{2(\sqrt{s+1}+1)}
\le\frac16.
\]
Therefore, the left-hand side of \eqref{eq:master-inequality} is bounded by
\begin{equation}\label{eq:G-def}
G(R):=
\frac{R-2}{2R}
+2\frac{R-2}{R-1}b(R)
+\frac{1+\Pi(R)}{\sqrt3}
+\frac16.
\end{equation}

Set
\begin{equation}\label{eq2}
q_0:=2+\frac1{\sqrt3},
\qquad
c_0:=2-\frac1{\sqrt3}.
\end{equation}
Since
\[
\Pi(R)=\frac{R}{R-1}b(R),
\]
we can rewrite \eqref{eq:G-def} as
\begin{equation}\label{eq:G-rewrite}
G(R)
=
\frac23+\frac1{\sqrt3}-\frac1R
+\left(q_0-\frac{c_0}{R-1}\right)b(R).
\end{equation}
Define
\begin{equation}\label{eq3}
F_3(R):=
\frac{17}{15\ln2}\ln R-G(R).
\end{equation}
Differentiating \eqref{eq:G-rewrite}, we obtain
\begin{equation}\label{eq:F3-derivative-exact}
\begin{aligned}
F_3'(R)
&=
\frac1R
\left[
\frac{17}{15\ln2}
-\frac{q_0-c_0/(R-1)}{\pi}
\right]
-\frac1{R^2}
-\frac{c_0\,b(R)}{(R-1)^2}.
\end{aligned}
\end{equation}

Consider
\[
f(R):=\frac{Rb(R)}{(R-1)^2}.
\]
Since
\[
f'(R)
=
\frac{R-1-(R+1)(\ln R+\eta)}
{\pi(R-1)^3}
<0,
\qquad R\ge4,
\]
the function $f$ is decreasing on $[4,\infty)$. Hence
\[
f(R)\le f(4),
\qquad R\ge4.
\]
Also,
\[
q_0-\frac{c_0}{R-1}\le q_0,
\qquad
\frac1R\le\frac14.
\]
Consequently, \eqref{eq:F3-derivative-exact} yields
\[
F_3'(R)
\ge
\frac1R
\left[
\frac{17}{15\ln2}
-\frac{q_0}{\pi}
-\frac14
-c_0f(4)
\right].
\]

Now
\begin{equation}\label{eq4}
b(4)=\frac{\ln4+\eta}{\pi}\approx0.481,
\end{equation}
and therefore
\[
f(4)
=
\frac{4b(4)}9
\approx0.214.
\]
Using \eqref{eq2}, we obtain
\[
\frac{17}{15\ln2}
-\frac{q_0}{\pi}
-\frac14
-c_0f(4)
\approx0.261>0.
\]
Thus,
\[
F_3'(R)>0,
\qquad R\ge4,
\]
and hence $F_3$ is increasing on $[4,\infty)$.

It remains to check $F_3(4)$. Since
\[
\frac{17}{15\ln2}\ln4
=
\frac{34}{15}
\approx2.267,
\]
and, by \eqref{eq:G-rewrite} and \eqref{eq4},
\[
\begin{aligned}
G(4)
&=
\frac23+\frac1{\sqrt3}-\frac14
+\left(q_0-\frac{c_0}{3}\right)b(4)\approx2.006,
\end{aligned}
\]
we obtain from \eqref{eq3}
\[
\begin{aligned}
F_3(4)
&=
\frac{34}{15}-G(4)\approx
2.267-2.006\approx0.2603
\end{aligned}
\]
and
$$F_3(4)>0.25>0.$$
Therefore,
$
F_3(R)\ge F_3(4)>0,
\ R\ge4.
$
This proves \eqref{eq:master-inequality} for $s\ge3$.

The three cases complete the proof of \eqref{eq:master-inequality}.
\end{proof}

\section{Proof of the main theorem}

\begin{proof}[Proof of Theorem~\ref{theorem1}]
We first assume $R\ge3$ and argue by induction on $n$.  For $n=1$,
$K_1\equiv1$, and
\[
1\le\frac{17}{15\ln2}\ln R.
\]
Assume the assertion is true for all positive integers smaller than $n$.
Choose $k$ such that $M_k\le n<M_{k+1}$ and write
$
n=sM_k+q,
\ 1\le s<m_k,
\ 0\le q<M_k.
$
By Lemmas~\ref{lemma:L1-reduction} and \ref{lemma:scalar},
\[
\begin{aligned}
n\|K_n\|_1
&\le q\|K_q\|_1
+sM_k\frac{17}{15\ln2}\ln R.
\end{aligned}
\]
If $q>0$, then $q<n$ and the induction hypothesis gives
\[
q\|K_q\|_1
\le q\frac{17}{15\ln2}\ln R;
\]
for $q=0$ this term is zero.  Since $n=sM_k+q$, we obtain
\[
\|K_n\|_1\le\frac{17}{15\ln2}\ln R.
\]
Thus the estimate holds for every $R\ge3$.

It remains to consider $R=2$.  Then $m_j=2$ for every $j\ge0$, so the
Vilenkin system is the Walsh--Paley system.  Toledo \cite{Toledo2018} proved
that
\begin{equation}\label{eq:Toledo}
\sup_{n\ge1}\|K_n\|_1=\frac{17}{15}.
\end{equation}
Since
\[
\frac{17}{15}=\frac{17}{15\ln2}\ln2,
\]
the asserted estimate also holds for $R=2$.  Hence the theorem recovers Toledo’s Walsh–Paley result in the endpoint case $R=2$ and extends the same sharp logarithmic estimate to all bounded Vilenkin groups.

Finally, if a universal constant $C$ satisfied
\[
\|K_n\|_1\le C\ln R
\]
for every bounded Vilenkin group and every $n\ge1$, then the Walsh--Paley
case and \eqref{eq:Toledo} would give
\[
\frac{17}{15}\le C\ln2.
\]
Thus $C\ge17/(15\ln2)$, proving optimality.
\end{proof}

\section*{Acknowledgments}

George Tephnadze was supported by Shota Rustaveli National Science Foundation (SRNSF), grant no. FR-24-698. The authors acknowledge the use of AI tools. All mathematical arguments and proofs in the final manuscript were checked and written by the authors.

\end{document}